\documentstyle[12pt]{article}
\def\be{\begin{equation}}
\def\ee{\end{equation}}
\def\lb{\label}
\begin{document}
\vskip 1.5cm
\begin{flushright}
{\bf Preprint IFUP-TH  62/95}
\end{flushright}
\vskip 1.5cm
\vspace{2cm}
\begin{center}
{\Large \bf Twisted Yang-Baxter equations for linear
quantum (super)groups}

\vspace{1.5cm}

{\bf A.P.Isaev}\footnote{Permanent address:
Bogoliubov Laboratory of Theoretical Physics, JINR,
Dubna, Moscow Region, 141 980, Russia }  \\
Dipartimento di Fisica, Universit\'{a} di Pisa, \\
Piazza Torricelli 2, 56100 Pisa, Italy

\end{center}

\vspace{2cm}

{\bf Abstract} \\

\vspace{0.5cm}
We consider the modified (or twisted) Yang-Baxter equations
for the $SL_{q}(N)$ groups and $SL_{q}(N|M)$ supergroups.
The general solutions for these equations are presented
in the case of the linear quantum (super)groups.
The introduction of spectral parameters in the twisted Yang-Baxter equation
and its solutions are also discussed.

\newpage

{\bf 1. Introduction} \\

Recently various types of modified Yang-Baxter equation (m-YBE)
have been considered. First of all such m-YBE
appeared in investigations on special exchange algebras \cite{CG}.
Another one was explored \cite{KS} in the context
of the construction of new integrable lattice models which generalize
the $SL(3)$ (and in general $SL(N)$)
chiral Potts model. After a series of
papers (see \cite{3D}) devoted to the
solutions of the tetrahedron equation,
similar modification, but now for the 3-d analogs of YBE,
have been used for a construction of
new integrable 3-d lattice theories \cite{3DD}.
Further, a variant
of YBE  has been found as certain cubic relations (for
$R$-matrices being a special set of quantum 6-j symbols)
which express the consistency of a quadratic algebra for
elements of matrix generating a set of Clebsh-Gordon coefficients
\cite{F}, \cite{AF},
\cite{BF}. It is interesting that this modification
coincides with the one considered in \cite{CG} and the corresponding
m-YBE and $R$-matrices essentially
depend on the phase space coordinates.
Note that the same dependence occures
for the classical $r$-matrices (which are called dynamical $r$-matrices)
in Calogero-Moser type
models (see \cite{Sk}, \cite{ABT} and references therein).
On the other hand, we recall that the q-analogs of the 6-j symbols
(or Racah coefficients) give the braiding/fusing matrices
expressing the property of crossing symmetry for
four-point conformal blocks
in 2-D conformal field theories (see e.g. \cite{AG}).
Finally we stress that the analogous `twisted' YBE has been
proposed also in the context of quasitriangular Hopf algebras \cite{D}.

In this paper we investigate the
m-YBE appeared in \cite{CG} as consistency relations for
exchange matrices and has been considered in \cite{F},
\cite{AF}, \cite{BF} as some relations for $SL_{q}(2)$ 6-j symbols.
Here, in the cases of $SL_{q}(N)$ and $SL_{q}(N|M)$ (super)groups,
we present the explicit solutions $R(p)$ for such m-YBE.
These solutions could be related with
6-j symbols for corresponding quantum groups. Then we show how one can
generalize these m-YBE and their solutions by introducing spectral
parameters and also present the
Yangian type limits for these solutions.
Our conjecture is that after introducing the spectral parameter
we obtain some objects related to the 6-j symbols
for quantum affine Kac-Moody algebras.

\vspace{1cm}
{\bf 2. Quantum deformations
of dynamical systems on co-adjoint
orbits or Alekseev-Faddeev toy models }
\vspace{0.5cm}

At the beginning,
to introduce the objects which will be under consideration,
we remind some facts from the paper \cite{AF}.
It is known \cite{OP} that apparently all
finite dimensional integrable models,
like Toda chain or Calogero particles, can be considered as systems
on co-adjoint orbits of some Lie groups $G$
(with Lie algebras ${\cal G}$) and described by the Lagrangians:
\be
\lb{1}
{\cal L}(t) = < L | \frac{d}{dt} g \, g^{-1} > - \frac{1}{2} <L|L> +
<L - \mu^{L}|\phi > +
<g^{-1}L g - \mu^{R}|\psi > \; ,
\ee
where $t$ is a time, $g(t) \in G$, $L(t)$
and constant elements $\mu^{L,R}$ belong to the space
${\cal G}^{*}$ dual to the Lie algebra ${\cal G}$, terms with
$\phi, \; \psi \in {\cal G}$ define the momentum mappings
and $\phi, \; \psi$
are nothing but Lagrange multipliers, \\
$<.|.>$ is a paring of
${\cal G}$ and ${\cal G}^{*}$.
We also identify ${\cal G}$ and ${\cal G}^{*}$ through
the invariant Killing metric. The explicit choice of the
group $G$, multipliers
$\phi$, $\psi$ and elements $\mu^{L,R}$ specifies the dynamical
system. If we consider the case for which we can take the matrix
representation for ${\cal G}$ and ${\cal G}^{*}$
such that the pairing will be defined via operator
$Tr$ $(<A|B> \rightarrow Tr(AB))$,
then, one can find the equations of motion from the Lagrangian (\ref{1})
and prove that the quantities $I_{n}=Tr(L^{n})$ are integrals of
motion. Thus, for appropriate momentum mappings we can expect that
the system with Lagrangian (\ref{1}) yields an example of integrable model.
 From the Lagrangian (\ref{1}) we find the following
Poisson brackets (see e.g. \cite{AF}):
\be
\lb{2}
\begin{array}{c}
\{ g^{1}, g^{2} \} =0 \; , \\
\{ L^{1}, \, g^{2} \} = C \, g^{2} \; , \\
\{ L^{1}, \, L^{2} \} = - \frac{1}{2} [ C, \, L^{1} - L^{2} ] \; .
\end{array}
\ee
where as usual $g^{1} = g \otimes 1, \; L^{2} = 1 \otimes L, \dots$ and
$C = t_{a} \otimes t_{b} \, \eta^{ab}$ is an ad-invariant tensor
($\eta^{ab}$ define Killing metric, and $t_{a}$ form the basis
for Lie algebra ${\cal G}$).
Then, for the case $G=SL(N)$, the diagonalization of the
left $L$ and right $g^{-1}Lg$ momenta can be considered:
\be
\lb{3}
L = u P u^{-1} \; , \;\; g^{-1}Lg = v^{-1} P v \; ,
\ee
and this leads to the diagonalization of the group element $g$:
\be
\lb{4}
g = u \, Q^{-1} \, v \; .
\ee
Here we have used
\be
\lb{5}
P=-\frac{i}{2}
{\rm diag}\{ p_{1}, \, p_{2}, \dots , p_{N} \} \; , \;\;
\sum_{i=1}^{N} p_{i} = 0 \; ,
\ee
\be
\lb{6}
Q={\rm diag}\{ \exp(i x_{1}), \, \exp(i x_{2}), \dots ,
\exp(i x_{N}) \} \; , \;\;
\sum_{i=1}^{N} x_{i} = 0
\ee
and matrices $u, \; v$ belong to the homogeneous space $G/H$ where
$H$ is a Cartan subgroup associated with $P$.

In the papers \cite{AF}, \cite{BF}, \cite{AT} it has been shown that Poisson
structure (\ref{2}), in terms of the new variables
$\{u, \, v, \, P, \, Q \}$, acquires the form
\be
\lb{7}
\{ u^{1}, \, u^{2} \} = - u^{1} \, u^{2} \, r_{0}(p) \; , \;\;
\{ v^{1}, \, v^{2} \} =  r_{0}(p) \, v^{1} \, v^{2} \; ,
\ee
and
\be
\lb{8}
\{ x_{i}, \, p_{j} \} = \delta_{ij} \; (1 \leq i,j \leq N-1)
\; , \;\; \{u_{0}, \, v_{0} \} =0 \; , \;\;
\{ u_{0}, \, p_{i} \} = 0 = \{ v_{0}, \, p_{i} \} \; ,
\ee
where we have introduced $u=u_{0}\, Q, \; v = Q \, v_{0}$,
\be
\lb{9}
r_{0}(p) = \sum_{\alpha} \frac{i}{p_{\alpha}}
(e_{\alpha} \otimes e_{-\alpha} - e_{-\alpha} \otimes e_{\alpha} )
\ee
and $\alpha$ runs over positive roots of ${\cal G}$. Namely we have
$$
e_{\alpha} = e_{jk}, \; j < k \; , \;\; (e_{ij}e_{kl} = \delta_{kj}e_{il})
\; , \;\; p_{\alpha} = (p_{j} - p_{k}) \; .
$$
The variable $p$ (in $r_{0}(p)$) means that $r_{0}$ depends on
all moments $p_{i}$.
The quantum version of the formulas
(\ref{7}) - (\ref{9}) has been discussed in \cite{AF}, \cite{BF} for the
case of $SL_{q}(2)$ group and, as it has been
pointed out in \cite{BF}, can be postulated for the general
case of $SL_{q}(N)$ in the same form:
\be
\lb{10}
\begin{array}{c}
R_{12} \, u_{1} \, u_{2} = u_{2} \, u_{1} \, R(p)_{12} \; , \\
R(p)_{12} \, v_{2} \, v_{1} = v_{1} \, v_{2} \, R_{12} \; ,
\end{array}
\ee
\be
\lb{10a}
\, [ u^{1}_{0}, \, v^{2}_{0} ] = 0 \; , \;\;
\, [ u_{0}, \, p_{i} ] = 0 = [ v_{0}, \, p_{i} ] \; , \;\;
[x_{i}, \, p_{j} ] = i \, h \, \delta_{ij} \; (i,j \leq N-1) \; .
\ee
Here $q$ - is a deformation parameter, $h$ is a Planck constant,
$R_{12}$ is the well known $R$-matrix for the $GL_{q}(N)$ group
(see \cite{J}, \cite{FRT}) and we introduce new $R$-matrix
$R(p)_{12}$ which nontrivially depends on the moments $p_{i}, \;
\forall i$. For simplicity we remove from eqs. (\ref{10})
nonessential factor $(q^{-1/N})$ which transforms $GL_{q}(N)$
$R$- matrix to the $SL_{q}(N)$ one.
Here and below we use $R$- matrix formalism which was developed in \cite{FRT}.
We note that $u-$ and $v-$ algebras (\ref{10}) can be identified
via relation $u = v^{-1}$.
Let us recall that the $GL_{q}(N)$ $R$ -matrix
satisfies the YBE:
\be
\lb{11a}
\hat{R} \, \hat{R}' \, \hat{R} = \hat{R}' \, \hat{R} \, \hat{R}'
\ee
and the Hecke relation:
\be
\lb{11}
\hat{R}_{12}^{2} = \lambda \hat{R}_{12} + 1 \; , \;\; \lambda =
q-q^{-1} \; .
\ee
where $\hat{R}= \hat{R}_{12} = P_{12} R_{12}$,
$\hat{R}' = P_{23} R_{23}$ and $P_{12}$ is a permutation matrix.
Using relation (\ref{11}) we immediately derive from eqs. (\ref{10}) that
$\hat{R}(p)=$ $\hat{R}(p)_{12}= P_{12} \, R(p)_{12}$ also obeys
the Hecke relation:
\be
\lb{12}
\hat{R}(p)^{2} = \lambda \hat{R}(p) + 1 \; ,
\ee
Considering third order monomials in $u$ (or in $v$)
and using the commutation relations (\ref{10}), (\ref{10a})
give the analogue of the YBE \cite{CG}, \cite{F}, \cite{AF}, \cite{BF}
for the new objects $R(p_{i})_{12}$:
\be
\lb{13}
\begin{array}{c}
(Q_{1})^{-1} \, R(p)_{23} \, Q_{1} \, R(p)_{13} \,
(Q_{3})^{-1} \, R(p)_{12} \, Q_{3} = \\
R(p)_{12} \, (Q_{2})^{-1} \, R(p)_{13} \, Q_{2} \, R(p)_{23}
\end{array}
\ee
This equation can be rewritten in the form
(cf. with (\ref{11a})):
\be
\lb{14}
\hat{R}(p) \, \widetilde{R}(p)' \, \hat{R}(p) =
\widetilde{R}(p)' \, \hat{R}(p) \, \widetilde{R}(p)'
\ee
where the matrix $\widetilde{R}(p)' = Q_{3} \, \hat{R}(p)_{23} \, (Q_{3})^{-1}$
obviously also satisfies the Hecke condition.
For searching solutions of equations (\ref{14}) it is convenient
to rewrite them as
\be
\lb{14a}
(Q^{-1}_{3} \, \hat{R}(p) \, Q_{3})
\, \hat{R}(p)' \,
(Q^{-1}_{3} \, \hat{R}(p) \, Q_{3})
 = \hat{R}(p)' \,
(Q^{-1}_{3} \, \hat{R}(p) \, Q_{3})
 \, \hat{R}(p)'
\ee
where $\hat{R}(p)' = \hat{R}(p)_{23}$. We call eqs. (\ref{13}) -
(\ref{14a}) modified or twisted Yang-Baxter equations. Note that
>from eq. (\ref{14}) one can obtain the relations which are similar
to the reflection equations:
$$
L(p) \, \widetilde{R}(p)' \, L(p) \, \widetilde{R}(p)' =
\widetilde{R}(p)' \, L(p) \, \widetilde{R}(p)' \, L(p) \; ,
$$
$$
\widetilde{L}(p) \, \hat{R}(p) \, \widetilde{L}(p) \, \hat{R}(p) =
\hat{R}(p) \, \widetilde{L}(p) \, \hat{R}(p) \, \widetilde{L}(p) \; ,
$$
where $L(p) = \hat{R}(p)^{2}, \; \widetilde{L}(p) = (\widetilde{R}'(p))^{2}$.

In the papers \cite{F}, \cite{AF} and \cite{BF}
the explicit form for the matrix $R(p)$ in $SL_{q}(2)$ case
has been presented.
There was stressed also that the elements of the
matrix $R(p)$ give the special set
of the 6-j symbols for $SL_{q}(2)$ group. Below we resolve equations
(\ref{14}) (and, thus, derive the
explicit formulas for $R(p)$) in the case of groups $SL_{q}(N)$ and
supergroups $SL_{q}(N|M)$.
Note, that special solutions of (\ref{14a}) for arbitrary
simple quantum groups (including $SL_{q}(N)$) are presented in
\cite{CG}.
Our solutions are multiparametric and more general.
The specific choices of these parameters lead
our solutions to the $R(p)$ which
could be interpreted as corresponding 6-j symbols \cite{F}
or as exchange $R$- matrix for the vertex algebra \cite{CG}.

\vspace{1cm}
{\bf 3. Solutions of the modified YBE for the case
of the linear quantum groups.}
\vspace{0.5cm}

The m-YBE and their solutions for the case of
$SL_{q}(N)$ have been firstly considered in \cite{CG}. Here we obtain
more general multiparametric solution for the case of linear
quantum groups and supergroups. We will consider the case of
$GL_{q}(N)$ and $GL_{q}(K|N-K)$.
In the case of
special q-groups, solutions $R(p)$ can be obtained by multiplying
$R(p)$ on some factor which is a simple function of $q$ (see below).

Let us search the solution of m-YBE (\ref{14a}) in the form
\be
\lb{15}
\hat{R}_{12} =
\hat{R}(p)^{i_{1}i_{2}}_{j_{1}j_{2}} =
\delta^{i_{1}}_{j_{2}}\delta^{i_{2}}_{j_{1}} \, a_{i_{1}i_{2}}(p) +
\delta^{i_{1}}_{j_{1}}\delta^{i_{2}}_{j_{2}} \, b_{i_{1}i_{2}} (p) \; .
\ee
Without limitation of generality one can put $b_{ii}(p) = 0$. Now the
condition that
$R(p)$ (\ref{15}) satisfy the Hecke relation
(\ref{12}) gives
the following constraints:
\be
\lb{16}
b_{ij} + b_{ji} = \lambda \; , \;\; i \neq j \; ,
\ee
\be
\lb{17}
a_{ij} \, a_{ji}  - b_{ij} \, b_{ji} = 1 , \; i \neq j \; ,
\ee
\be
\lb{18}
a_{i}^{2} - \lambda a_{i} - 1 = 0 \Rightarrow
a_{i} - \lambda = \frac{1}{a_{i}}
\; , \;\; a_{i} \equiv a_{ii} \; .
\ee
Note that the eq. (\ref{18}) has two solutions:
$a_{i} = \pm q^{\pm 1}$ and therefore coefficients $a_{i}$ are
independent of the parameters $p_{k}$. If we take
$a_{i} = q \; , \; \forall i$
(or $a_{i} = -q^{-1} \; , \; \forall i$) then we will have the case
of the group $GL_{q}(N)$ (or $GL_{-q^{-1}}(N)$). But if we consider
the mixing case: $a_{i} = q$ for $ 1 \leq i \leq K$ and
$a_{i} = - q^{-1}$ for $ K + 1 \leq i \leq N$ then we come to the
case of supergroups $GL_{q}(K|N-K)$.

Now let us use the relations (cf. with (\ref{10a}))
\be
\lb{19a}
\exp(-i \, x_{j}) \, p_{k} \, \exp(i \, x_{j}) = p_{k} + h \delta_{kj} \; ,
\;\; (1 \leq k,j \leq N)
\ee
and substitute
(\ref{15}) to the m-YBE (\ref{14a}). As a result, in addition to the relations
(\ref{16}) - (\ref{18}), we obtain new constraints on  the functions
$a_{ij}(p)$ and $b_{ij}(p)$. First of all we deduce:
\be
\lb{19}
a_{ij}(p_{1}, \dots , p_{N}) = a_{ij}(p_{i},p_{j}) \; , \;\;
b_{ij}(p_{1}, \dots , p_{N}) = b_{ij}(p_{i},p_{j}) \; , \;\;
\ee
and, as a consequence, relations (\ref{16}), (\ref{17})
have to be rewritten in the form:
$$
b_{ij}(p_{i},p_{j}) + b_{ji}(p_{j},p_{i}) = \lambda \; ,
$$
$$
a_{ij}(p_{i},p_{j}) \, a_{ji}(p_{j},p_{i}) -
b_{ij}(p_{i},p_{j}) \, b_{ji}(p_{j},p_{i}) = 1 \; .
$$
Then we have the constraint:
\be
\lb{20}
b_{ij} \, b_{jk} \, b_{ki} +
b_{ik} \, b_{kj} \, b_{ji} =0 \; , \;\; i \neq j \neq k \neq i
\ee
(where there is no summation over $i,j,k$) and equations
\be
\lb{21}
b_{ij}(p_{i} + h, p_{j}) = \frac{ b_{ij}(p_{i},p_{j}) \, a_{i}}{1/a_{i} +
b_{ij}(p_{i},p_{j})} \; ,
\ee
\be
\lb{22}
b_{ij}(p_{i} , p_{j}+ h) = \frac{  b_{ij}(p_{i},p_{j}) / a_{j}}{
a_{j} - b_{ij}(p_{i},p_{j})} \; .
\ee
Using (\ref{21}) and (\ref{22}) leads to the following
general relations:
$$
b_{ij}(p_{i} + n \, h, p_{j} + m \, h) =
\frac{a_{i}^{n} \, a_{j}^{-m} \, b_{ij}(p_{i},p_{j}) }{
a_{i}^{-n}\, a_{j}^{m} +
 b_{ij}(p_{i},p_{j}) \, (a_{i}^{n}a_{j}^{-m}
-a_{i}^{-n}a_{j}^{m})/\lambda } =
$$
\be
\lb{23}
= \frac{ \lambda \, a_{i}^{n} \, a_{j}^{-m} \, b_{ij}(p_{i},p_{j}) }{
a_{i}^{n}\, a_{j}^{-m}\, b_{ij}(p_{i},p_{j}) +
a_{i}^{-n} \, a_{j}^{m} \, b_{ji}(p_{j},p_{i})}
\; .
\ee
$\,$ From these relations one can immediately find the general
solution for the coefficients $b_{ij}(p)$:
\be
\lb{24}
b_{ij}(p_{i},p_{j}) =
\frac{a_{i}^{p_{i}/h} \, a_{j}^{-p_{j}/h} \,
b_{ij}^{0} }{a_{i}^{-p_{i}/h}\, a_{j}^{p_{j}/h}
+ b_{ij}^{0} \, (a_{i}^{p_{i}/h}a_{j}^{-p_{j}/h} -
a_{i}^{-p_{i}/h}a_{j}^{p_{j}/h})/ \lambda } =
\ee
$$
= \frac{ \lambda \, a_{i}^{p_{i}/h} \, a_{j}^{-p_{j}/h}
\, b_{ij}^{0} }{
a_{i}^{p_{i}/h}\, a_{j}^{-p_{j}/h}\, b_{ij}^{0} +
a_{i}^{-p_{i}/h} \, a_{j}^{p_{j}/h} \, b_{ji}^{0} } \; ,
$$
where constants $b_{ij}^{0} = b_{ij}(0,0)$ have to
obey the algebraical relations:
\be
\lb{25}
\begin{array}{c}
b^{0}_{ii}=0 \; , \;\; b^{0}_{ij} + b^{0}_{ji} = \lambda \; ,  \\ \\
b^{0}_{ij} \, b^{0}_{jk} \, b^{0}_{ki} +
b^{0}_{ik} \, b^{0}_{kj} \, b^{0}_{ji} =0 \; ,
\end{array}
\ee
which can be deduced by the substitution (\ref{24}) into
eqs. (\ref{16}) and (\ref{20}).

It is clear now that if
$a_{i} = a_{j}$ (indices $i$ and $j$ `have the same grading')
then $b_{ij}(p_{i},p_{j}) = b_{ij}(p_{i}-p_{j})$, but if
$a_{i} = -1/a_{j}$ (the case of supergroups when indices
$i$ and $j$ `have the opposite grading') then we deduce that
$b_{ij}(p_{i},p_{j}) = b_{ij}(p_{i}+p_{j})$.
Note that the only conditions on the parameters $a_{ij}(p)$
needed for fulfillment of m-YBE
are listed in (\ref{17}) and (\ref{19}).

Let us consider the case of the $GL_{q}(N)$ group. In this
case we have $a_{i} = q \; \forall i$ and relation (\ref{24})
takes the form:
\be
\lb{26}
b_{ij}(p_{i}-p_{j}) =
\frac{q^{(p_{i} -p_{j})/h} \, b_{ij}^{0} }{
q^{(-p_{i} + p_{j})/h} + [ \frac{p_{i}-p_{j}}{h} ]_{q} \, b_{ij}^{0} } \; ,
\ee
while for the functions $a_{ij}(p_{i} , p_{j})$ we obtain
>from (\ref{17}) the following relations:
\be
\lb{27}
a_{ij}(p_{i} , p_{j}) \, a_{ji}(p_{j} , p_{i}) =
1 + \frac{\lambda^{2} \, b_{ij}^{0} \, b_{ji}^{0}}{
(q^{(p_{i} - p_{j})/h} \, b_{ij}^{0} +
q^{(p_{j} - p_{i})/h} \, b_{ji}^{0} )^{2} } \; .
\ee
In (\ref{26}) we have used the standard notation
$[x]_{q} = (q^{x} - q^{-x})/ (q-q^{-1})$.
One can obtain from these expressions the solutions discussed in \cite{CG},
\cite{F}, \cite{AF} and \cite{BF} if we
take the proper normalization of $a_{ij}$ (e.g. Faddeev's
or unitary normalization
$a_{ij}(p_{i}-p_{j}) = a_{ji}(p_{j}-p_{i})$) and
consider, in (\ref{26}), (\ref{27}),
the limit $b_{ij}^{0} \rightarrow \infty$ $(i<j)$. This limit can be
performed selfconsistently such that it does not cancel the
conditions (\ref{25}).
Then we recall that our consideration was done for the
case of the general groups $GL_{q}(K|N-K)$. In fact one can
obtain $SL$-matrices $R(p)$ multiplying them by the functions
$q^{1/(N-K) - 1/K}$ which are needed for obtaining the
identity $Sdet_{q}(R) = 1$ (for the $SL_{q}(N)$ case
we have to multiply $R(p)$ by $q^{-1/N}$).
As a result for the $SL_{q}(N)$ case we have
the matrix $[q^{-1/N} \cdot R(p)]$ (\ref{15}) with substitution:
\be
\lb{*}
b_{ij}(p_{i}-p_{j}) = \frac{q^{(p_{i}-p_{j})/h} }{
[(p_{i}-p_{j})/h]_{q}} = \lambda - b_{ji}(p_{j}-p_{i}) \; ,
\ee
\be
\lb{**}
a_{ij}(p_{i}-p_{j}) = \frac{ \left(
[ (p_{i}-p_{j})/h+1 ]_{q} \, [ (p_{i}-p_{j})/h - 1 ]_{q} \right)^{1/2} }{
[(p_{i}-p_{j})/h]_{q} \, \epsilon_{ij} } = a_{ji}(p_{j}-p_{i})
\ee
where $\epsilon_{ij} = \pm 1 \; (i < j), \; = \mp 1 \; ( j > i)$.
This choice of $a_{ij}$ leads to the unitary condition (for real
$q$ and $p_{i}^{\dagger} = p_{i}$):
$$
R(p)_{12}^{\dagger} = R(p)_{12}^{t_{1}t_{2}} = R(p)_{21}.
$$

Analogously, for the $SL_{q}(K|N-K)$ case, we obtain the matrix \\
$q^{1/(N-K) - 1/K} \cdot R(p)$ (\ref{15}):
\be
\lb{?}
\begin{array}{l}
\hat{R}(p)_{SL_{q}(K|N-K)} = q^{1/(N-K) - 1/K}
\left(
\delta^{i_{1}}_{j_{2}} \, \delta^{i_{2}}_{j_{1}}
\left[
(-1)^{(i_{1})} \, q^{1-2(i_{1})} \delta^{i_{1}i_{2}} +
\right.
\right.
\\ \\
\left.
 + a(p_{i} - p_{j}) ( \theta_{K+1,i} \, \theta_{K+1,j} +
\theta_{i,K}\theta_{j,K} ) +
a(p_{i} + p_{j}) ( \theta_{K+1,i} \, \theta_{j,K} +
\theta_{i,K}\theta_{K+1,j} )
\right] +
\\ \\
+ \delta^{i_{1}}_{j_{1}} \, \delta^{i_{2}}_{j_{2}}
\left[
b(p_{i} - p_{j})  \theta_{K+1,i} \, \theta_{K+1,j} +
b(p_{j} - p_{i})  \theta_{i,K} \, \theta_{j,K} +
\right.
\\ \\
\left.
\left.
b(p_{i} + p_{j})  \theta_{K+1,i} \, \theta_{j,K} +
b(- p_{i} - p_{j})  \theta_{i,K} \, \theta_{K+1,j}
\right] \;
\right)
\end{array}
\ee
where $(i) = 0$ for $1 \leq i \leq K$ and $=1$
for $K+1 \leq i \leq N$, $\theta_{ij} = 1$ for $i>j$ and $=0$
for $i \leq j$. The functions $a(p_{i}-p_{j}) = a_{ij}(p_{i}-p_{j})$,
$b(p_{i}-p_{j}) = b_{ij}(p_{i}-p_{j})$ are defined in
(\ref{*}), (\ref{**}).

To conclude this section we stress that we have found the more general
solution of the m-YBE (\ref{13}) - (\ref{14a})
then obtained in \cite{CG}, \cite{F}.
Namely, our solutions depend on the
set of arbitrary parameters $b_{ij}^{0}$ constrained by the conditions
(\ref{25}). The role of these parameters is still to be clarified.
Taking some special limits and choosing the normalization of $a_{ij}$ one
leads to the known solutions $R(p)$ of the papers \cite{CG} and \cite{F}
(see e.g. (\ref{*}), (\ref{**})).

\vspace{1cm}
{\bf 4. Modified (twisted) YBE with spectral parameters.}
\vspace{0.5cm}

In this section we show that m-YBE (\ref{13}) - (\ref{14a})
can be generalized by introducing spectral
parameters $y, \, z, \dots$. We demonstrate that every
solution $R(p)$ which has been found in the previous section will lead to
the solution $R(p,y)$ for the m-YBE with the spectral parameters.
It is interesting to note that such introducing of the
spectral parameters  can be done in
complete analogy with the usual way of obtaining
the trigonometric solutions
\be
\lb{tr}
\hat{R}(y) = y^{-1} \, \hat{R} - y \, \hat{R}^{-1}
\ee
of the YBE from the
$R$- matrices $\hat{R}$ related to the $GL_{q}(N)$ groups.
On the other hand, we know that the trigonometric solutions (\ref{tr})
are related to the quantum Kac-Moody algebras \cite{Dr}.
In this connection (following the statements of the papers
\cite{F}, \cite{AF}, \cite{BF})
it is natural to conjecture that $R(p,\, y)$ could be interpreted as
a special set of 6-j symbols for the q-deformations
of linear affine algebras.

The natural assumption about the form of the m-YBE
dependent on the spectral parameters
is the following:
\be
\lb{28}
\hat{R}(p,\, y) \, \widetilde{R}'(p,\, y \cdot z) \, \hat{R}(p,\, z)
=\widetilde{R}'(p,\, z) \, \hat{R}(p,\, y \cdot z) \, \widetilde{R}'(p,\, y)
\ee
Now it is not difficult to check by using m-YBE (\ref{14})
and the Hecke relations (\ref{12}) that the following
matrices
$$
\hat{R}(p,\, y) = y^{-1} \, \hat{R}(p) - y \, \hat{R}(p)^{-1}
\; , \;\;
$$
\be
\lb{29}
\widetilde{R}'(p,\, y) =
y^{-1} \, \widetilde{R}'(p) -
y \, (\widetilde{R}'(p))^{-1} =
Q_{3} \, \hat{R}'(p, \, y) \, Q_{3}^{-1} \; , \;\;
\ee
are the solutions of the new m-YBE (\ref{28}).
We note that solutions (\ref{29}) satisfy the identity
$$
\hat{R}(p, \, y) \, \hat{R}(p, \, y^{-1}) =
\left( \lambda^{2} - (y-y^{-1})^{2} \right) \; ,
$$
which is a kind of unitary condition
for $R(p)$ (if $y^{*} = y^{-1}$). It is clear that the analog
of the relations (\ref{10}) dependent on the spectral
parameters has the form
$$
\hat{R}(y) \, u_{1}(yz) \, u_{2}(z) =
u_{1}(z) \, u_{2}(yz) \, \hat{R}(p, \, y) \; .
$$
Now let us put $q=exp(\gamma \, h)$ \cite{AF}, \cite{BF} and
$y=exp((-1/2)\lambda \, \theta)$.
Following \cite {AF} we consider
two different cases: deformed classical
case ($h = 0, \; \gamma \neq 0$) and
quantum non-deformed case ($\gamma =0, \; h \neq 0$). In the first
case we obtain that $\hat{R}(p, \, y)/\lambda$ tends to the
Yangian $R$- matrix $\hat{R}(\theta)=\theta \, P_{12} -1$
which satisfies usual YBE
$$
\hat{R}(\theta) \, \hat{R}'(\theta + \theta') \, \hat{R}(\theta')
=\hat{R}'(\theta') \, \hat{R}(\theta + \theta') \, \hat{R}'(\theta) \; .
$$
In the second case we derive:
\be
\lb{29a}
\lim_{\gamma \rightarrow 0} \frac{\hat{R}(p, \, y)}{\lambda}
= \hat{R}(p, \, \theta)
 = \theta \, \hat{R}^{0}(p)  - 1 \; ,
\ee
where $\hat{R}^{0}(p)$ is represented in the form (\ref{15})
with the following parameters (cf. with \cite{AF,BF}):
$$
b_{ij} = \frac{h}{p_{i}-p_{j}} = - b_{ji} \; ,
$$
$$
a_{ij} = \frac{((p_{i}-p_{j} +h)(p_{i}-p_{j} -h))^{1/2}}{
\epsilon_{ij}(p_{i} - p_{j})} = a_{ji} \; ,
$$
and matrix $\hat{R}(p, \, \theta)$ (\ref{29a}) satisfies twisted
YBE:
\be
\lb{29b}
\hat{R}(p,\, \theta) \,
\widetilde{R}'(p,\, \theta + \theta') \, \hat{R}(p,\, \theta') =
\widetilde{R}'(p,\, \theta') \, \hat{R}(p,\, \theta + \theta') \,
\widetilde{R}'(p,\, \theta) \; .
\ee

To conclude this paper we note that it would be
extremely interesting to use
the twisted YBE (\ref{28}), (\ref{29b}) and
their solutions (\ref{29}), (\ref{29a}) for formulating the
integrable models e.g. via box construction of \cite{KS}
or to relate these solutions with the braiding matrices
describing generalized statistics (see e.g. \cite{LM}).

\vspace{1cm}
{\bf Acknowledgments}
\vspace{0.5cm}

The author would like to thank R.Kashaev for helpful comments and
 M.Mintchev for valuable discussions and
the idea of introducing spectral parameter in the twisted YBE.
I am grateful to A.Di Giacomo and M.Mintchev
for hospitality at the Physics Department of the University
of Pisa where this
paper was done. This work was partially supported by
the exchange program between INFN and JINR (Dubna), ISF (grant
RFF 300) and RFFI (grant 95-02-05679-a).

\end{document}